\documentstyle[amstex,amssymb]{amsart}

\newcommand{\Zint}{{\mathbb {Z}}}    
     
\newcommand{\Rea}{{\mathbb {R}}}      
     
\newcommand{\Ha}{{\frak {h}}}    
\newcommand{\halmos}{\rule{5pt}{5pt}}

\numberwithin{equation}{section}

\newtheorem{prop}{\bf Proposition}[section]
\newtheorem{thm}[prop]{\bf Theorem}
\newtheorem{lemma}[prop]{\bf Lemma}
\newtheorem{cor}[prop]{\bf Corollary}
\newtheorem{con}{\bf Conjecture}
\newenvironment{rmk}{\noindent{\bf Remark}\hskip 5pt}{\hfill{$\Box$}}

\begin{document}

\title[Eigenstates of the Calogero-Moser model]
{On the eigenstates of the elliptic Calogero-Moser model}
\author{Kouichi Takemura}
\address{Research Institute for Mathematical Sciences, Kyoto University, 606-8502 Kyoto, Japan.}
\email{takemura@@kurims.kyoto-u.ac.jp}
\thanks{K.T. is supported by JSPS Research Fellowship for Young Scientists} 
\dedicatory{To the memory of Denis Uglov}


\begin{abstract}
It is known that the trigonometric Calogero-Sutherland model is obtained by the trigonometric limit $(\tau \rightarrow \sqrt{-1} \infty)$ of the elliptic Calogero-Moser model, where $(1,\tau)$ is a basic period of the elliptic function.

We show that for all square-integrable eigenstates and eigenvalues of the Hamiltonian of the Calogero-Sutherland model, if $\exp (2\pi \sqrt{-1} \tau )$ is small enough then there exist square-integrable eigenstates and eigenvalues of the Hamiltonian of the elliptic Calogero-Moser model which converge to the ones of the Calogero-Sutherland model for the $2$-particle and the coupling constant $l$ is positive integer cases and the $3$-particle and $l=1$ case.

In other words, we justify the regular perturbation with respect to the parameter $\exp (2\pi \sqrt{-1} \tau )$.

With some assumptions, we show analogous results for $N$-particle and $l$ is positive integer cases. 
\end{abstract}

\maketitle

\section{Introduction}

The elliptic Calogero-Moser model (or elliptic Olshanetsky-Perelomov model) is a quantum many body system whose Hamiltonian is given as follows,
\begin{equation}
 H :=- \frac{1}{2} \sum_{i=1}^{N} 
\frac{\partial ^{2}}{\partial x_{i}^{2}}
+ l (l+1)
\sum_{1 \leq i<j \leq N}
\wp ( x_{i}-x_{j}),
\end{equation}
where $\wp(x)$ is the Weierstrass elliptic function. (\cite{OP})

This model is known to be integrable, i.e. there exists $N$-algebraically independent commuting operators which commute with the Hamiltonian $H$.  (\cite{OP})

In this article, we are going to investigate the eigenstates of the Hamiltonian $H$.

It is known that the Hamiltonian of the trigonometric Calogero-Sutherland model is obtained by the trigonometric limit $(\tau \rightarrow \sqrt{-1} \infty)$ of the elliptic Calogero-Moser model, where $(1,\tau)$ is a basic period of the elliptic function.

For the trigonometric case, the eigenstates are known. They are described by the Jack polynomial. The eigenvalues are also known.

The main idea of this article is to connect the elliptic model with the trigonometric model,
and obtain some information about the eigenstates and the eigenvalues of the elliptic model. 

Since the first term in the expansion of the Hamiltonian of the elliptic model as a power series in $p=\exp (2\pi \sqrt{-1} \tau )$ is up to constant the Hamiltonian of the trigonometric model, we can obtain the formal eigenstates of the elliptic model by the perturbation method, i.e. we obtain the formal eigenstates as the formal power series of $p$.
But the convergence is not guaranteed a priori.
For example the formal perturbation expansion of the eigenstates of $\tilde{H} =-\frac{d^2}{dx^2}+x^2+\alpha x^4$ with respect to $\alpha $ does not convergent regularly.

The main result of this article gives a sufficient condition for the regular convergence of the perturbation expansion.
In particular, for the $2$-particle and the coupling constant $l$ is positive integer cases and the $3$-particle and $l=1$ case, we have the convergence for all eigenstates related to Jack polynomial.
The theorem guarantee sthe numerous square-integrable eigenstates and their eigenvalues. (see Conjecture \ref{maincon}, Remark below Conjecture \ref{maincon}, and Theorem \ref{mth1})

The main tool is the Bethe Ansatz method.
The Bethe Ansatz method replaces the problem of finding eigenfunction of the Hamiltonian by solving the transcendental equation which is called the Bethe Ansatz equation.
For the $2$-particle and $l \in \Zint_{>0} $ cases, this reduces to the Bethe Ansatz for the Lam\'e equation and was performed more than a century ago. (\cite{WW})
For the $N$-particle and  $l \in \Zint_{>0} $ cases, Felder and Varchenko obtained the Bethe Ansatz equation by investigating the asymptotic behavior of the integral representation of the solution of the Knizhnik-Zamolodchikov-Bernard equation. (\cite{FVKZB})

After obtaining the Bethe Ansatz equation, there are two things to be considered.
The first one is to find the condition when the eigenfunction obtained by the Bethe Ansatz method is connected to the square-integrable eigenstate and the second one is how the solution of the Bethe Ansatz equation behaves.

The first question is not so difficult. The condition is described as a certain continuous parameter belonging to some lattice. For details see Lemma \ref{silemma}.

The second question is serious.
We consider the solution at $p=0$ (the trigonometric limit) and look into the behavior where $p$ is near $0$.
In this step, the key lemma is the implicit function theorem.

In this way, we construct the square-integrable eigenstates and obtain the main result.

Let us comment on the relationship between the eigenstate obtained by the Bethe Ansatz method and the finite dimensional invariant subspace preserved by the Hamiltonian of the elliptic Calogero-Moser model.

For the case $l\in \Zint _{>0}/N$, the Ruijsenaars operators (a $q$-analogue of the operators of the Calogero-Moser model) preserve the finite dimensional subspace of theta-type function, which depends on $l$. 
By considering the limit $q\rightarrow 1$, we recognize that the commuting operators of the elliptic Calogero-Moser model preserve the finite dimensional subspace of periodic functions, which depends on $l$.
(See \cite{BFV,Che,Has,HK} etc.)

The space spanned by the square-integrable eigenstate obtained by the Bethe Ansatz method is different from the the finite dimensional space of doubly periodic functions.
If a function obtained by the Bethe Ansatz method belong to the finite dimensional space of doubly periodic functions, the function has poles and is not square-integrable.
To obtain the square-integrable eigenstates, we will consider the (anti-)symmetrization in section \ref{trigl}. In the procedure of (anti-)symmetrization, the function vanishes.

Since the eigenstate obtained by the Bethe Ansatz method is directly connected to the Jack (or Macdonald) polynomial, we can conclude that the diagonalization of the finite dimensional space of the theta-type function is not directly connected to the Jack (or Macdonald) polynomial with a ``physical'' parameter but with a ``non-physical'' parameter.

This article is organized as follows.

In section 2, we review the properties of the Jack polynomial.
In section 3, we discuss the Bethe Ansatz method of the elliptic Calogero-Moser model, trigonometric limit and corresponding results.
In section 4, we solve the Bethe Ansatz equation, which is necessary to establish the main theorem. 
In section 5, we give some comments.

\section{Jack polynomial}

\subsection{Calogero-Sutherland model and Jack polynomial}

The Calogero Sutherland model is a model of a 1-dimensional quantum many body system.(\cite{Sut}) The Hamiltonian is given by
\begin{equation}
 H_{CS} :=- \frac{1}{2} \sum_{i=1}^{N} 
\frac{\partial ^{2}}{\partial x_{i}^{2}}
+\pi ^{2}
l(l+1)
\sum_{1 \leq i<j \leq N}
\frac{1}{\sin ^{2} (\pi  (x_{i}-x_{j}))}.
\end{equation}

We set $X_{i}:=\exp \left( 2 \pi \sqrt{-1} x_{i} \right)$ and ${\cal{M}}_N:=\{ \lambda = (\lambda _{1} \: , \: \lambda _{2} \: , \dots , \:  \lambda _{N}   )| i>j \Rightarrow \lambda _{i} - \lambda_j \in \Zint_{\geq 0} \}$.

The eigenstates of the Calogero-Sutherland model are described by the Jack polynomial $J_{\lambda}^{(\frac{1}{l+1})}(X)$ (\cite{Sta,Mac}), i.e.
\begin{equation}
 H_{CS} ( J_{\lambda}^{(\frac{1}{l+1})}(X)\Delta (X)^{l+1} ) = (e_0+ 2\pi^2 E_{\lambda }^{[\frac{1}{l+1}]})J_{\lambda}^{(\frac{1}{l+1})}(X)\Delta (X)^{l+1},
\end{equation}
where $e_0:= \frac{1}{6} \pi ^{2} (l+1)^{2} N(N^{2}-1)$, $\Delta (X):=(X_{1}X_{2} \dots X_{N} )^{\frac{1-N}{2}} \prod _{i<j} (X_{i}-X_{j})$ and $E_{\lambda }^{[\alpha ]}:= \sum_{i=1}^{N} \lambda _i^2 + \sum_{i=1}^{N} \frac{N+1-2i}{\alpha } \lambda _i $.
In particular, the ground-state is given by $\Delta (X)^{l+1}$. 

We set $|\lambda |:= \sum_{i=1}^{N}\lambda _i$ and define a partial ordering in ${\cal{M}}_N$ by 
$ \lambda \succeq \mu$ $\Leftrightarrow $ 
 $ |\lambda|=|\mu|, \; \sum_{j=1}^{i} \lambda_{j} \geq \sum_{j=1}^{i} \mu_{j}$ 
$ (i=1, \dots ,N)$. 

Let $m^{\lambda}$ be a monomial symmetric function associated with $\lambda$, i.e. $m^{\lambda}:= \sum_{\mu \in S_m\cdot \lambda}X_1^{\mu _1}\dots X_N^{\mu _N}$. The function $m^{\lambda}$ is a polynomial up to the multiplication of $(X_1 \dots X_N)^a$ for some $a$.

We summarize the property of the polynomial $J_{\lambda}^{(\alpha )}$.
\begin{equation}
J_{\lambda}^{(\alpha )}(X)=m^{\lambda }+ 
\sum_{\mu \prec \lambda } \tilde{c} ^{(\alpha )}_{\lambda , \mu } m^{\mu },
\label{Jnorma}
\end{equation}
\begin{equation}
\langle J_{\lambda}^{(\alpha )}(X), J_{\mu}^{(\alpha )} (X)\rangle = \delta_{\lambda, \mu}c_{\lambda},
\end{equation}
where the inner product $\langle \cdot , \cdot \rangle$ is given by
\begin{equation}  
\langle f, g \rangle := \frac{1}{N!}\left( \prod_{i=1}^N \oint_{|X_i|=1} 
\frac{dX_i}
{2\pi \sqrt{-1}X_i}  \right) \overline{\Delta (X)^{\frac{1}{\alpha}}f(X_1,\dots,X_N)} \Delta (X) ^{\frac{1}{\alpha}}g(X_1,\dots,X_N),
\end{equation}
, $\overline{X_i}= X_i^{-1}$ and $\tilde{c} ^{(\alpha )}_{\lambda , \mu }$, $c_{\lambda}$ are some constants.

We will define the Jack polynomial associated with the weight $\xi \in P^+$.
We put $\xi =\sum_{i=1}^{N} \xi_i \epsilon _i$ then $\xi_i -\xi_j \in \Zint _{>0}$ $(i>j)$, and we set $J_{\xi }^{(\alpha )}(X):=J_{(\xi_1,\dots ,\xi_N)}^{(\alpha )}(X)$.
The eigenvalue $E_{\lambda }^{[\alpha ]}$ is written as $(\lambda +\frac{1}{\alpha }\rho , \lambda +\frac{1}{\alpha }\rho )- \frac{1}{\alpha ^2}(\rho, \rho )$, where $\rho := \sum_{i=1}^{N} \frac{2N+1-2i}{2}\epsilon _i$ is the half sum of positive roots.

\section{Bethe Ansatz for the elliptic Calogero-Moser model}

\subsection{Bethe Ansatz}

In this section, we introduce the Bethe Ansatz method for the $N$-particle elliptic Calogero-Moser model with the coupling constant $l$ positive integer.
Most of the results mentioned in this section are due to Felder and Varchenko.(\cite{FVKZB, FVthr})

From this section we adopt the Hamiltonian shifted by some constant, i.e.
\begin{equation}
H^{\tau ,(l)} := -\frac{1}{2}\sum_{i=1}^{N} \frac{\partial^2}{\partial x_i^2}+
l(l+1)\sum_{1\leq i < j\leq N} (\wp(x_i -x _j)+2\eta ),
\end{equation}
where
$\eta := \pi^2 (\frac{1}{6}-4\sum_{n=1}^{\infty} \frac{p ^{n}}{1-p ^{n}})$ and $p=\exp (2\pi \sqrt{-1} \tau )$.

First we fix the parameters $N$ and $l$. We set $m:=lN(N-1)/2$.
Let $c: \{1,\dots ,m\} \rightarrow \{1, \dots ,N\}$ be the unique non-decreasing function such that $c^{-1}(j)$ has $(N-j)l$ elements.

Let $\epsilon _i$ $(1\leq i \leq N)$ be an orthonormal basis of $\Rea ^{N}$.
We realize the simple roots of $A_{N-1}$ type as $\alpha _i =\epsilon_i - \epsilon_{i+1}$.
We also realize the set of roots $R$, the root lattice $Q$, the weight lattice $P$, and the set of the dominant weights $P^+$ in the space $\Ha ^{*}= \{ \sum_{i=1}^{N}x _i \epsilon_i | \sum_{i=1}^{N}x _i=0 \}$

We set $p_i:=i(2N-i-1)l/2$ and define 
\[
V_{i}:= \{ p_{i-1}+1, p_{i-1}+2, \dots , p_i \} \; \; (1\leq i\leq N-1).
\] 
Let $W$ be the set of maps $w=(w_1, \dots ,w_N)$ $(w_i: \; V_i \rightarrow \{i,i+1,\dots ,N-1 \} )$ such that $\# \{ w_i^{-1}(j)\} =l$ for $1\leq i\leq j\leq N-1$.
For $w=(w_1, \dots ,w_{N-1}) \in W$, let $F_{w}$ be the set of maps $f=(f_1, \dots ,f_{N-2})$ $(f_{i}: \; V_{i+1} \rightarrow V_i)$ such that (i) $f_{i}$ is injective (ii) If $f_{i}(x)=y$ then $w_{i+1}(x)=w_i(y)$.

We set 
\begin{equation}
\theta_1(x):=2\sum_{n=1}^{\infty} (-1)^{n-1} \exp (\tau \pi \sqrt{-1}(n-1/2)^2) \sin(2n-1)\pi x, \label{th1}
\end{equation}
\begin{equation}
\theta (x):=\frac{\theta_1(x)}{\theta_1 '(0)}, \; \; \; \sigma_{\lambda }(x):=\frac{\theta '(0)\theta(x-\lambda)}{\theta(x)\theta(\lambda)}.
\end{equation}

We introduce the functions $\Phi _{\tau}(t_1, \dots , t_m)$ and $\omega (t;x)$ as follows
\begin{align}
& \Phi _{\tau}(t_1, \dots , t_m):= 
 e^{2\pi \sqrt{-1} (\xi,\sum_i t_j\alpha_{c(j)}) } \\
& \prod_{1\leq j \leq (N-1)l} \theta(t_j) ^{-lN} 
\prod_{c(i)=c(j)\atop{i< j}}\theta(t_i-t_j )^{2}
\prod_{|c(i)-c(j)|=1\atop{i< j}}\theta(t_i-t_j )^{-1}, \\
& \omega (t;x)=e^{2\pi \sqrt{-1} (\xi ,\sum_{i} x_i \epsilon_i)}
\sum_{w \in W} \sum_{f \in F_w} \prod_{i=1}^{N-1}  \prod_{k= p_{i-1}+1}^{p_i}
\sigma_{x_i -x_{w_i(k)+1}}(t_k-t_{f_i(k)}) \label{Bvect} 
\end{align}
where $t_0=0$, $f_0(k)=0$.

Then we have
\begin{prop} $($\cite{FVKZB, FVthr}$)$
If $(t^0_1, \dots , t^0_m)$ satisfy the following Bethe Ansatz equations,
\begin{equation}
\frac{\partial \Phi_{\tau}}{\partial t_i}| _{(t^0_1, \dots , t^0_m)}=0 \; \; (1\leq i \leq m).
\end{equation}
the function $\omega (t^0;x)$ is an eigenfunction of the Hamiltonian $H^{\tau ,(l)}$ with the eigenvalue 
\begin{equation}
2\pi ^2(\xi,\xi) -2\pi \sqrt{-1} \frac{\partial }{\partial \tau} S(t_{\tau, 1}^0 ,\dots,t_{\tau, N}^0 ; \tau ),
\label{CSev}
\end{equation}
where $S(t_1,\dots,t_N; \tau )= \sum_{i<j}(\alpha_{c(i)}, \alpha_{c(j)}) \log \theta (t_i -t_j)- \sum_{c(i)=1} lN\log \theta (t_i)$. 
\end{prop}
\begin{rmk}
At a glance, the expression of $\omega (t^0;x)$ is different from Felder and Varchenkos' one, but the two are the same in fact.
Our expression is indicated by the limit $(q\rightarrow 1)$ of the expression of the Bethe vector of the Ruijsenaars model, which is obtained by Billey.(\cite{Bil})
\end{rmk}
Therefore if we find the solution of the Bethe Ansatz equation, we can investigate the Calogero-Moser model more effectively.
For this purpose, we will consider the trigonometric limit and investigate the solutions of the Bethe Ansatz equations for both the elliptic and the trigonometric model.

\subsection{Trigonometric limit and the main theorem} \label{trigl}

We investigate the trigonometric limit to have a link to the Calogero-Sutherland model.

We set $T_i :=e^{-2\pi \sqrt{-1}t_i}$ and $X_i:=e^{2\pi \sqrt{-1}x_i}$.
As $p \rightarrow 0 $, we have the limits $\theta(x) \rightarrow \sin \pi x$, $\Phi_{\tau} (t_1, \dots , t_m)\rightarrow \mbox{const.}\Phi_{tri}(T_1, \dots , T_m)$, and $\omega(t,x) \rightarrow  \mbox{const.}\omega_{tri}(T,X)$, where
\begin{align}
& \Phi_{tri}(T_1, \dots , T_m) :=
\prod_{j=1}^{m} T_j^{-(\xi - \rho, \alpha_{c(j)})) }
\prod_{c(j)=1}(1-T_j)^{-lN }  \label{triP} \\
& \; \; \; \; \; \; \prod_{c(i)=c(j)\atop{i< j}}(T_i-T_j )^{2}
\prod_{|c(i)-c(j)|=1\atop{i< j}}(T_i-T_j )^{-1},\nonumber
\end{align}
\begin{equation}
\omega_{tri}(T,X)=
\frac{\prod_{i=1}^{N} X_i^{(\xi, \epsilon_i)}}
{\prod_{i<j }(X_i- X_j)^l}
\sum_{w \in W} \sum_{f \in F_w} \prod_{i=1}^{N-1}  \prod_{k= p_{i-1}+1}^{p_i}
\frac{X_i T_k-X_{w_i(k)+1}T_{f(k)}}{ T_k-T_{f(k)}},
\end{equation}
where $T_0 =1$ and $f_0 (k)=0$.

Let $F^{\tau}_{N,l}$ (resp. $F_{N,l}$) be a complement set of the set of zeros and poles of the function of the function
$\Phi_{\tau} (t_1, \dots , t_m)$ (resp. $\Phi_{tri}(T_1, \dots , T_m)$).

We set 
\[
\mbox{Sym}^{(l)}f(x_1, \dots ,x_N):= \left\{ 
\begin{array}{ll}
\sum_{\sigma \in S_N}f(x_{\sigma(1)}, \dots ,x_{\sigma(N)}) & l\mbox{ is odd,} \\
\sum_{\sigma \in S_N}\mbox{sgn}(\sigma) f(x_{\sigma(1)}, \dots ,x_{\sigma(N)}) & l\mbox{ is even.} 
\end{array}
\right. 
\]

The following conjecture describes the behavior of the solution of the trigonometric Bethe Ansatz equation. For some special cases, the conjecture is true.
\begin{con} \label{maincon}
Let $\tilde{\xi } = \sum_{i=1}^{N-1} m_i \Lambda_i = \sum_{i=1}^{N} m'_{i}\epsilon _i $ be an element of $P^+$. If $(\sum_{i=1}^{N-1} m_i \Lambda_i , \alpha ) \not \in  \{ 0,\pm 1,\dots ,\pm l \}$ for all $\alpha \in R$, there exists $\sigma \in S_N$ such that there is a non-degenerate critical point $(T^0_1, \dots, T^0_m) \in F_{N,l}$ of $\Phi_{tri}(T_1, \dots ,T_m)$ with the parameter $\xi = \sum_{i=1}^{N} m'_{\sigma (i)}\epsilon _i $ such that $\mbox{Sym}^{(l)}\omega_{tri}(T^0, X)$ is non--zero.
Here the non-degenerate critical point means the critical point $(T_1^0, \dots , T_m^0 )$(i.e. $\frac{\partial \Phi}{\partial T_i}|_{(T_1^0, \dots , T_m^0 )}=0)$ such that the Hessian at the critical point is non-zero.
\end{con}

\begin{rmk}
For each $N\in \Zint_{\geq 2}$ and $l\in \Zint_{>0}$, there exists infinite $\tilde{\xi } \in P^+$ such that Conjecture \ref{maincon} is proved.
\end{rmk}

\begin{prop} \label{prconj}
Conjecture \ref{maincon} is proved for the $N=2$ and $l \in \Zint_{>0}$ cases, and the $N=3$ and $l=1$ case.
\end{prop}
We will discuss this in detail for the $N=2$ and $l \in \Zint_{>0}$ cases in section \ref{a1} and he $N=3$ and $l=1$ case in section \ref{a2l1}

The following lemma is the key point to connect the trigonometric solutions and the elliptic ones.
\begin{lemma} \label{ift}
Let $\xi $ be an element of $\Ha ^*$. Let $(T_1^0,T_2^0, \dots ,T_m^0)\!$ $=(e^{-2\pi \sqrt{-1}t_1^0},e^{-2\pi \sqrt{-1}t_2^0},$ $\! \dots ,e^{-2\pi \sqrt{-1}t_m^0}) \in F_{N,l}$ be a non-degenerate critical point of $\Phi _{tri}$ (\ref{triP}). Then there exists some $\epsilon \in \Rea _{>0}$ such that if $|p |< \epsilon$ then there is  a non-degenerate critical point $(t_{\tau ,1}^0, t_{\tau ,2}^0,\dots , t_{\tau ,m}^0)$ of $\Phi _{\tau}$ (i.e. a non-degenerate solution of the elliptic Bethe Ansatz equation) for the same $\xi $ and $(t_{\tau ,1}^0, t_{\tau ,2}^0,\dots , t_{\tau ,m}^0)$ $\rightarrow $  $(t_1^0,t_2^0, \dots ,t_m^0)$ as $p \rightarrow 0$.
\end{lemma}
This lemma follows from the implicit function theorem. The condition $(T_1^0,T_2^0, \dots ,T_m^0) \in F_{N,l}$ is necessary to apply the implicit function theorem.

We introduce lemmas which are needed to obtain the main theorem.

\begin{lemma} \label{silemma}
Let $\xi = \sum_{i=1}^{N} m_i \Lambda_i $ be an element of the weight lattice $P$. Assume $(\sum_{i=1}^{N-1} m_i \Lambda_i , \alpha )\not \in  \{ 0,\pm 1,\dots ,\pm l \}$ for all $\alpha \in R$. Let $(t_{\tau,1}^0,t_{\tau,2}^0, \dots ,t_{\tau,m}^0) \in F^{\tau}_{N,l}$ be a solution of the Bethe Ansatz equation for $\Phi_{\tau}$ and $\omega(t^0_{\tau};x)$ be the Bethe vector (\ref{Bvect}). Then $\mbox{Sym}^{(l)}\omega (t^0_{\tau};x)$ is square-integrable on $[0,1]^N$ and also the eigenfunction of the operator $H^{\tau , (l)}$. 
\end{lemma}

\begin{rmk}
Lemma \ref{silemma} determines the condition for the existence of the square-integrable eigenstates for  $H^{\tau , (l)}$.
If $\xi = \sum_{i=1}^{N} m_i \Lambda_i \in (\Ha ^* \setminus P)$ then 
$\mbox{Sym}^{(l)}\omega (t^0_{\tau};x)$  is not square-integrable.
\end{rmk}

\begin{lemma} \label{BJ}
Let $\xi = \sum_{i=1}^{N} m_i \Lambda_i$ be an element of $P$. Assume $(\sum_{i=1}^{N-1} m_i \Lambda_i , \alpha ) \not \in  \{ 0,\pm 1,\dots ,\pm l \}$ for all $\alpha \in R$. 
Choose $\xi '\in P^+$ and $\sigma \in S_N$ such that $\xi' = \sigma (\xi)$.
 Let $(T_1^0,T_2^0, \dots ,T_m^0)$ be a solution of the Bethe Ansatz equation associated with $\xi$ and $\omega_{tri}(T^0;X)$ be the corresponding eigenfunction.
Then $\mbox{Sym}^{(l)}\omega _{tri}(T^0;X)= \mbox{const.} J^{(\frac{1}{l+1})}_{\xi '-(l+1)\rho }(X)\cdot \Delta (X)^{l+1}$.
\end{lemma}

By applying the Lemmas \ref{ift}, \ref{silemma}, \ref{BJ}, we have

\begin{thm} \label{mth1}
Let $\lambda $ be an element of $P^{+}$. We set $p =e^{2\pi \sqrt{-1} \tau }$. Assuming Conjecture \ref{maincon} for the weight $\lambda +(l+1)\rho $, then there exists $\epsilon \in \Rea _{>0}$ such that if $|p |< \epsilon$ then for each $p$ there exist eigenvalues $E^{\tau}_{\lambda}$ and eigenfunctions $F_{\lambda }^{\tau}(x_1,\dots ,x_N)$ of the Hamiltonian of the elliptic Calogero-Moser model  $H^{\tau, (l)}$ such that  
$E^{\tau}_{\lambda} \rightarrow 2 \pi ^{2}E_{\lambda}^{[\frac{1}{l+1}]}+e_0 + \frac{\pi^2}{6}N(N-1)l(l+1)$, $F_{\lambda }^{\tau}(x_1, \dots x_N) \rightarrow 
\Delta (X) ^{l+1} J^{(\frac{1}{l+1})}_{\lambda }(X_1, \dots ,X_N)$ as $p \rightarrow 0$ and $F_{\lambda }^{\tau}(x_1,\dots ,x_N) $ is square-integrable, i.e.
\[
\int_0^1 \cdots \int_0^1 |F^{\tau}_{\lambda }(x_1,\dots x_N)|^2 dx_1\dots dx_N < \infty.
\]
Here the function $J^{(\frac{1}{l+1})}_{\lambda }(X_1, \dots ,X_N)$ is the Jack polynomial associated to the weight $\lambda$ and $X_i=\exp (2\pi \sqrt{-1} x_i)$.
\end{thm}

As a corollary, we have
\begin{cor} \label{mcor1}
For the $N=2$ and $l \in \Zint_{>0}$ cases, and the $N=3$ and $l=1$ case, Theorem \ref{mth1} is proved for all $\lambda \in P^+$ without assumption.
\end{cor}

From the uniqueness of the perturbation expansion up to constant, we have
\begin{thm}
Let $\lambda $ be an element of $P^{+}$.
If the condition of Theorem \ref{mth1} holds for $\lambda $,
the perturbation expansion related to the eigenstate $J^{(\frac{1}{l+1})}_{\lambda }(X_1, \dots ,X_N)$ converges regularly.

In particular, for the $N=2$ and $l \in \Zint_{>0}$ cases, and the $N=3$ and $l=1$ case, the perturbation expansion related to the eigenstate $J^{(\frac{1}{l+1})}_{\lambda }(X_1, \dots ,X_N)$ converges regularly for all $\lambda \in P^{+}$
\end{thm}

\section{Solutions of the Bethe Ansatz equation}
\subsection{$2$-particle case} \label{a1}

In this section we will consider the $N=2$ case.

The Bethe Ansatz equation is given by $\frac{\partial \Phi_{\tau}}{\partial t_i}=0$ $(1\leq i \leq l)$, where
\[
\Phi _{\tau}=
 e^{\pi \sqrt{-1}\sum_{j=1}^{l} m_1 t_l  } \prod_{1\leq j \leq l} \theta(t_j) ^{-2l} 
\prod_{1\leq i< j\leq l}\theta(t_i-t_j )^{2}.
\]

We set $s:=x_1-x_2$. The operator $H^{\tau} $ is replaced by
\begin{equation}
H^{\tau} = -\frac{d^2}{ds^2}+l(l+1)(\wp(s)+2\eta). \label{CMH2}
\end{equation}
We calculate the l.h.s. of the equation (\ref{Bvect}). We find that the eigenfunction (Bethe vector) is equal to
\begin{equation}
e^{\pi \sqrt{-1} m_1 s}\frac{\theta (s-t_1) \dots \theta (s-t_l)}{\theta (s)^l},
\label{Bv2}
\end{equation}
up to constant.

We will confirm Conjecture \ref{maincon} for the $N=2$ and $l\in \Zint_{>0}$ case.
For this purpose, we will investigate the trigonometric case.
In this case the condition $(\sum_{i=1}^{N-1} m_i \Lambda_i , \alpha ) \not \in  \{ 0,\pm 1,\dots ,\pm l \}$ for $\forall \alpha \in R$ is equivalent to the condition $m_1 \not \in \{ 0,\pm 1,\dots ,\pm l \}$.

The function $\Phi_{tri}$ for the $N=2$ case is 
\[
\Phi _{tri}=
\prod_{j=1}^{l} T_j^{(-m_1+1) }
(1-T_j)^{-2l }
\prod_{1\leq i<j \leq l}(T_i-T_j )^{2}.
\]

The critical points of the function $\Phi_{tri}$ and the Hessian at the critical points were calculated by Varchenko to calculate the norms of the Bethe vector which is obtained by the asymptotic behavior of the integral representation of the solution of the KZ equation.

Let $(T_1^0 , \dots ,T_l^0)$ be a critical point of $\Phi $. Set 
\[ 
\sigma _i := \sum_{j_1< \dots <j_i} T_{j_1}^0 \dots T_{j_i}^0, \; \; \; \;
\overline{\sigma }_i :=\sum_{j_1< \dots <j_i} (1-T_{j_1}^0) \dots (1-T_{j_i}^0).
\]
Then we have
\begin{prop}$($\cite{V}$)$
\[
\sigma_i = 
\left(
\begin{array}{c}
l \\
i
\end{array}
\right)
\prod_{j=1}^{i}
\frac{-m_1+1+l-j}{-m_1-j},
\; \; \; \; 
\overline{\sigma }_i = 
\left(
\begin{array}{c}
l \\
i
\end{array}
\right)
\prod_{j=1}^{i}
\frac{l+j}{m_1+j}.
\]
Moreover the critical point exists uniquely if the numerators and the denominators are non-zero for all $i \in \{1, \dots ,l \}$.
\end{prop}

We set
\[
\delta := \prod_{1\leq i<j \leq l}(T_i^0-T_j^0)^2, \; \;
\mbox{Hess}:= 
\left.
\mbox{det}
\left(
\left(
\frac{\partial ^2}{\partial {T_i}\partial {T_j}}(-\kappa \log \Phi)
\right)
_{i,j}
\right)
\right|
_{T=T^0}
\]

\begin{prop}$($\cite{V}$)$ \label{discHess}
\end{prop}
\[
\delta = \prod_{j=0}^{l-1}
\frac{(j+1)^{j+1} (-m_1+1+j)^j (-2l+j)^j}{(-m_1-j-1)^{2l-j-2}},
\]
\[
\mbox{Hess}=l!\prod_{j=0}^{l-1}
\frac{(-m_1-j-1)^3}{(-m_1+1+j)(-2l+j)}.
\]


If $m_1 \not \in \{ 0,\pm 1,\dots ,\pm l \}$ then $\sigma_i, \overline{\sigma }_i (i=1,\dots ,l), \delta, \mbox{Hess}$ are finite and non-zero. Hence the solutions of the Bethe Ansatz equation belong to $F_{2,l}$ and are non-degenerate.

The function $\omega_{tri}(T,X)$ is given by
\begin{equation}
\omega_{tri}(T,X)
= \mbox{const.}\frac{ X_1^{\frac{1}{2}m_1} X_2^{-\frac{1}{2}m_1}}
{(X_1- X_2)^l}
\prod_{k=1}^{l}
(X_1 T_k -X_2 ).
\end{equation}

We can check that if $m_1 \not \in \{ 0,\pm 1,\dots ,\pm l \}$ then $\mbox{Sym}^{(l)}\omega_{tri}(T,X)$ is non-zero.

Therefore we have proved Conjecture \ref{maincon} for the $N=2$ and $l \in \Zint_{>0}$ case and obtain Corollary \ref{mcor1}.

\subsection{$3$-particle and $l=1$ case} \label{a2l1}

We will consider the $N=3$ and $l=1$ case.

The Bethe Ansatz equation is given by $\frac{\partial \Phi_{\tau}}{\partial t_i}=0$ $(1\leq i \leq 3)$, where
\[
\Phi _{\tau}=
 e^{2\pi \sqrt{-1}( m_1 t_1+m_1t_2+m_2t_3 )  }  ((\theta(t_1) \theta(t_2))^{-3}
\theta(t_1-t_2 )^{2}
(\theta(t_1-t_3 )\theta(t_2-t_3 ))^{-1}.
\]

The eigenstate (the Bethe vector) of the operator $H^{\tau , (l)}$ with $N=2$ and $l=1$ is equal to
\begin{align}
& e^{2\pi \sqrt{-1} (\frac{m_1}{3}(2x_1-x_2-x_3)+\frac{m_2}{3}(x_1+x_2-2x_3))}
\nonumber \\
& \left\{
\frac{\theta(x_1-x_2-t_1)}{\theta(t_1)}
\frac{\theta(x_1-x_3-t_2)}{\theta(t_2)}
\frac{\theta(x_2-x_3-t_3+t_2)}{\theta(t_3-t_2)} \right.
+
\nonumber \\
& \left.
\frac{\theta(x_1-x_3-t_1)}{\theta(t_1)}
\frac{\theta(x_1-x_2-t_2)}{\theta(t_2)}
\frac{\theta(x_2-x_3-t_3+t_1)}{\theta(t_3-t_1)}
\right\}
,
\nonumber
\end{align}
up to constant.

We will confirm Conjecture \ref{maincon} for the $N=3$ and $l=1$ case.
For this purpose, we will investigate the trigonometric case.
In this case the condition $(\sum_{i=1}^{N-1} m_i \Lambda_i , \alpha ) \not \in  \{ 0,\pm 1,\dots ,\pm l \}$ for $\forall \alpha \in R$ is equivalent to the condition $m_1 ,m_2,m_1+m_2\not \in \{ 0,\pm 1 \}$.

The function $\Phi_{tri}$ is 
\[
\Phi_{tri} =
(T_1T_2)^{(-m_1+1) }
(T_3)^{(-m_2+1) }
((1-T_1)(1-T_2))^{-3 }
(T_1-T_2)^{2}
(T_1-T_3)^{-1}
(T_2-T_3)^{-1}
\]

The critical point of $\Phi_{tri} $ is 
\[
T_3^0=\frac{(m_1+m_2-1)(m_2-1)}{(m_1+m_2+1)(m_2+1)},
\]
and $(T_1^0,T_2^0)$ are the solution of the following equation,
\[
(m_1+m_2+1)(m_1+1)X^2+2(-m_1^2-m_1m_2+2)X+(m_1+m_2-1)(m_1-1)=0.
\]
Therefore we have
\[
T_1^0T_2^0 = \frac{(m_1+m_2-1)(m_1-1)}{(m_1+m_2+1)(m_1+1)},
\]
\[
(1-T_1^0 )(1-T_2^0 )= \frac{6}{(m_1+m_2+1)(m_1+1)},
\]
\[
\prod_{1\leq i<j\leq 3}(T_i^0-T_j^0)^2=
\frac{2(m_1+m_2-1)^2(2m_1^2+2m_1m_2-m_2^2-3)^3}
{(m_1+1)^4(m_2+1)^4(m_1+m_2+1)^6}.
\]

The value of the Hessian at the critical point is
\[
\frac{1}{6}
\frac{(m_1+1)^3(m_2+1)^3(m_1+m_2+1)^5}{(m_1-1)(m_2-1)(m_1+m_2-1)^3},
\]

Hence if $m_1 , m_2, m_1 +m_2\not \in \{ 0,\pm 1 \}$, then the solutions of the Bethe Ansatz equation with the parameter $(m_1, m_2)=( \tilde{m_1} , \tilde{m_2} )$ or  $(m_1, m_2)=(-\tilde{m_2},-\tilde{m_1})$ belong to $F_{3,1}$ and are non-degenerate.

The function 
$\omega_{tri}(T,X)$ is given by
\begin{align}
& \omega_{tri}(T,X)
= \mbox{const.}\frac{ X_1^{\frac{2\tilde{m}_1}{3}+\frac{\tilde{m}_2}{3}} X_2^{\frac{-\tilde{m}_1}{3}+\frac{\tilde{m}_2}{3}} X_3^{\frac{-\tilde{m}_1}{3}-\frac{2\tilde{m}_2}{3}}}
{(X_1- X_2)(X_2-X_3)(X_2-X_3)} \nonumber \\
& 
\left\{
\frac{(X_1T_1-X_2)}{T_1}
\frac{(X_1T_2-X_3)}{T_2}
\frac{(X_1T_3-X_2T_2)}{T_3-T_2} 
\right.
 \nonumber \\
& 
\left.
 +
\frac{(X_1T_1-X_3)}{T_1}
\frac{(X_1T_2-X_2)}{T_2}
\frac{(X_1T_3-X_2T_1)}{T_3-T_1}
\right\}
.\nonumber
\end{align}

We can also check that if $m_1 , m_2, m_1 +m_2\not \in \{ 0,\pm 1 \}$ and  $(m_1, m_2) \in F_{3,1}$ then $\mbox{Sym}^{(l)}\omega_{tri}(T,X)$ is non-zero.

Therefore we have proved Conjecture \ref{maincon} for the $N=3$ and $l=1$ case and obtain Corollary \ref{mcor1}.

\section{Concluding remarks}
\subsection{}
In this article we have explained how to obtain the square-integrable eigenstates based on the Bethe Ansatz method and have shown the convergence of the perturbation series.

We will comment on some features of the perturbation expansion.

The elliptic Hamiltonian admits the expansion $H=H_0+\sum_{i=1}^{\infty}p^iV_i$, where $H_0$ is the trigonometric Hamiltonian up to constant.
To get the perturbation expansion, we must calculate $V_i\psi _{\lambda}= \sum c^{(i)}_{\lambda, \mu}\psi_{\mu}$, where $\psi _{\lambda}$ are the eigenstates of $H_0$.
In our cases, this process is reduced to the calculation of the Pieri formula of the Jack polynomial, i.e. $J_{\lambda} e_j = \sum d^{(i)}_{\lambda, \mu}J_{\mu}$, where $e_j$ is the $j$-th elementary symmetric polynomial.
The Pieri formula is written explicitly in \cite{Mac}.

There are merits for the perturbation method. The calculation of the perturbation does not essentially depend on the coupling constant $l$ although the calculation of the Bethe Ansatz method strongly depends on $l$.
The Bethe Ansatz method is valid for $l \in \Zint_{>0}$ but the perturbation method is valid for all cases if we ignore the convergence.
We hope that the ``eigenstates'' obtained by the perturbation method converge for all cases.

\subsection{}
The q-deformation of the elliptic Calogero Moser model is known.
It is called the elliptic Ruijsenaars model.
For the elliptic Ruijsenaars model, the Bethe Ansatz method is also known. (\cite{FVeqg,Bil})

Based on the Bethe Ansatz method, we can obtain similar results on the elliptic Ruijsenaars model.

\vspace{.15in}
{\bf Acknowledgment}
The author would like to thank Professors M. Kashiwara and T. Miwa for discussions and support. Thanks are also due to Dr. Tim Baker.

\end{document}